\documentclass{article}
\usepackage[english]{babel}
\usepackage{amsmath}
\usepackage{amssymb}
\usepackage{amsthm}
\usepackage{url}
\usepackage{amsfonts}

\usepackage[all]{xy}

\DeclareMathOperator{\id}{Id}
\DeclareMathOperator{\im}{Im} 
\DeclareMathOperator{\kr}{Ker} 
\DeclareMathOperator{\cil}{Cylinder} 
\newcommand{\W}{{B}} 
\DeclareMathOperator{\barc}{Bar} 

\newcommand{\rrdc}{\mbox{\,\(\Rightarrow\hspace{-7pt}\Rightarrow\)\,}}
\newcommand{\lrdc}{\mbox{\,\(\Leftarrow\hspace{-7pt}\Leftarrow\)\,}}
\newcommand{\lrrdc}{\mbox{\,\(\Leftarrow\hspace{-7pt}\Leftarrow\hspace{-5pt}\Rightarrow\hspace{-7pt}\Rightarrow\)\,}} 

\newcommand{\Zset}{\mathbb{Z}}         
\newcommand{\Nset}{\mathbb{N}}

\newtheorem{thm}{Theorem}

\newtheorem{alg}{Algorithm}

\theoremstyle{definition}
\newtheorem{defn}{Definition}
\newtheorem{rmk}{Remark}

\title{Computing the homology of groups: \\the geometric way \footnote{This research was partly supported by Ministerio de Educaci\'on y Ciencia (Spain), project MTM2009-13842-C02-01.}}

\author{A. Romero and J. Rubio}

\date{Departamento de Matem\'aticas y Computaci\'on\\
 Universidad de La Rioja \\
 Edificio Vives, c/Luis de Ulloa s/n \\
 26004 Logro\~no, Spain. \\
 ana.romero@unirioja.es \\
 julio.rubio@unirioja.es  }

\begin{document}

\maketitle

\begin{abstract}

In this paper we present several algorithms related with
the computation of the homology of groups, from a geometric
perspective (that is to say, carrying out the calculations
by means of simplicial sets and using techniques of
Algebraic Topology). More concretely, we have developed some algorithms which, making use
of the \emph{effective homology} method, construct
the homology groups of Eilenberg-MacLane
spaces $K(G,1)$ for different groups $G$, allowing one in particular to determine the
homology groups of $G$.

Our algorithms have been programmed as new modules for
the Kenzo system, enhancing it with the following new
functionalities:
\begin{itemize}
  \item construction of the effective homology of $K(G,1)$ from
  a given finite free resolution of the group $G$;
  \item construction of the effective homology of $K(A,1)$ for every
  finitely generated Abelian group $A$ (as a consequence, the effective
  homology of $K(A,n)$ is also available in Kenzo, for all $n \in \mathbb{N}$);
  \item computation of homology groups of some $2$-types;
  \item construction of the effective homology for central extensions.
\end{itemize}
In addition, an \emph{inverse} problem is also approached in this work:
given a group $G$ such that $K(G,1)$ has effective homology, can a finite free
resolution of the group $G$ be obtained? We provide some algorithms to solve
this problem, based on a notion of \emph{norm} of a group, allowing us to
control the convergence of the process when building such a resolution.
\end{abstract}

\section{Introduction}

When using homological algebra techniques to study group theory, two
different (but related) alternatives are possible (see \cite{Bro82}
for details on the following discussion). One is algebraic
and is based on the notion of \emph{resolution} (replacing the group
under study by an acyclic object in a suitable category of modules).
The other alternative is geometric in nature. It consists in finding
a contractible topological space with a free action of a group
$G$. Then the space of orbits of the action can be endowed with a
convenient quotient topology, in such a way that we obtain an
\emph{aspherical} space (that is to say, a space whose only non-null
homotopy group is the first, fundamental one). The homology of this
space is, by definition, the homology of $G$, and it does not depend
on the choosing of the contractible space or of the action. Each aspherical
space (unique up to homotopy type) is a particular \emph{Eilenberg-MacLane
space} for $G$, and is generically denoted by $K(G,1)$.

If we move to \emph{computational} mathematics, the preferred via chosen
was the algebraic one, as exemplified by the package HAP
(\cite{HAP}) of the computer algebra system GAP (\cite{GAP}), which contains
an impressive number of algorithms dealing with resolutions. The geometric
way has been up to now neglected from the algorithmic point of view.
The reason is that the contractible spaces to be constructed are very
frequently of infinite type (even in cases where the group $G$ is not
too complicated), apparently closing the possibility of a computational
treatment.

This view changed drastically when Sergeraert introduced at the
end of the 1980's his theory of \emph{effective homology} (\cite{Ser94}).
His methods allow the programmer to deal with spaces of infinite dimension,
encoded in a lazy functional programming style, producing a complete revision
of Algebraic Topology from a constructive point of view (see \cite{RS06} for
recent developments of this theory). Perhaps more important from a practical
point of view was Sergeraert's construction of the \emph{Kenzo} system, a
Common Lisp program implementing the effective homology methods (\cite{Kenzo}).
Since then, the programmer can work on a computer with simplicial sets,
loop spaces, fibrations, classifying spaces and many other Algebraic Topology
constructions, computing, at the end, homology groups of complicated spaces
(under the combinatorial form of \emph{simplicial sets}).

Taking into account this new situation, this paper represents a first step
to take up again the geometrical way of approaching \emph{group homology}, by means of
techniques from effective homology and using, and extending, Kenzo as a
computing platform.

Our proposal is not opposed to the algebraic view. Our aim is rather to take
the best of both worlds. Therefore, and as a first module, we programmed,
in collaboration with Graham Ellis (see \cite{RER09}), an OpenMath link
between HAP and Kenzo, allowing Kenzo to import from HAP resolutions of groups.
Once a resolution of a group $G$ is internally stored in Kenzo, an algorithm
allows us to construct the Eilenberg-MacLane space $K(G,1)$, with
\emph{effective homology}. This provides not only access to some homology
groups of $G$, but also makes it possible to apply on the space $K(G,1)$ all the
powerful tools available in Kenzo, and construct in this way further spaces.

This via is explored in this paper. We show two applications in Algebraic
Topology, and another one in Homological Algebra. As a first application, we
develop a Kenzo package to compute, as \emph{objects with effective homology},
the generalized Eilenberg-MacLane spaces $K(G,n)$ for any finitely generated
Abelian group $G$ and for all $n \in \mathbb{N}$. These objects are very
important in Algebraic Topology, to study and compute homotopy groups,
through Whitehead and Postnikov towers (see \cite{May67} and \cite{RS06}).

As a second topological application, we compute mechanically (for the first
time, up to our knowledge) some homology groups of $2$-types, the second
step (the first one consists of Eilenberg-MacLane spaces) towards the difficult
problem of characterizing homotopy types.

Our last application provides programs to deal with the effective homology of
central extensions of groups. The theoretical algorithms were known some time
ago (see \cite{Rub97}), but only now the technological tools explained before
allow us to tackle the problem of programming them. Let us observe that this \emph{algebraic}
application has also positive consequences on topological problems, since it
enlarges the field of application of our $2$-types package: we can also compute
with 2-types whose fundamental group is a central extension.

Finally, we have also approached an \emph{inverse} problem: how to obtain a
resolution of a group $G$ from the knowledge of an effective homology of
$K(G,1)$. The results in this area are still partial, and more research will
be needed to get fully satisfactory algorithms, and to proceed to implement them
as Kenzo modules.

The organization of the paper is as follows. The next section is devoted to
preliminaries. Section 3 contains our main algorithm, which constructs the effective
homology of $K(G,1)$ from a finite resolution of a group $G$, and then Section 4
collects some interesting fields of application of this result. Section 5
explains how our algorithms have been translated to Common Lisp and comments on
experimental results. In Section 6 an inverse problem is considered: given a group
$G$ with effective homology, it is (sometimes) possible to determine a resolution for
$G$. The paper ends with conclusions, open problems and the
bibliography.

\section{Definitions and preliminaries}

\subsection{Some fundamental notions about homology of groups}
\label{sec:group-hmlg}

The following definitions and important results about homology of groups can be found in~\cite{Mac63} and~\cite{Bro82}.

\begin{defn}
\label{defn:chcm}
Given a ring $R$, a \emph{chain complex} of \mbox{$R$-modules} is a pair of sequences \mbox{$C_\ast=(C_n,d_n)_{n \in \Zset}$} where, for each degree $n\in \Zset$, $C_n$ is an $R$-module and $d_n:C_n \rightarrow C_{n-1}$ (\emph{the differential map}) is an $R$-module morphism such that \mbox{$d_{n-1}\circ d_n= 0$} for all $n$.
\end{defn}

\begin{defn}
Let $C_\ast=(C_n,d_n)_{n \in \Zset}$ be a chain complex of $R$-modules, with $R$ a general ring. For each degree \mbox{$n\in \Zset$}, the \emph{$n$th homology module} of $C_\ast$ is defined to be the quotient module
$H_n(C_\ast)=\kr d_n / \im d_{n+1}$.
A chain complex $C_\ast$ is \emph{acyclic} if $H_n(C_\ast)=0$ for all~$n $.
\end{defn}


\begin{defn}
Let $G$ be a group and $\Zset G$ the free $\Zset$-module generated by the elements of $G$. The multiplication in $G$ extends uniquely to a $\Zset$-bilinear product $\Zset G\times \Zset G\rightarrow \Zset G$ which makes $\Zset G$ a ring. This is called the \emph{integral group ring} of $G$.
\end{defn}

\begin{defn}
A \emph{resolution} $F_\ast$ for a group $G$ is an acyclic chain complex of $\Zset G$-modules
$$
\cdots \longrightarrow F_2 \stackrel{d_{2}}{\longrightarrow} F_1 \stackrel{d_1}{\longrightarrow} F_0  \stackrel{\varepsilon}{\longrightarrow} F_{-1}=\Zset \longrightarrow 0
$$
where $F_{-1}=\Zset$ is considered a $\Zset G$-module with the trivial action and $F_i=0$ for $i<-1$. The map $\varepsilon: F_0 \rightarrow F_{-1}=\Zset$ is called \emph{augmentation}.
If $F_i$ is free for each $i\geq 0$, then $F_\ast$ is said to be a \emph{free resolution}.

Very frequently, resolutions come equipped with a \emph{contracting homotopy} $h$, which is a set of Abelian group morphisms $h_n:F_n \rightarrow F_{n+1}$ for each $n \geq -1$ (in general not compatible with the $G$-action), such that
\begin{equation*}
\begin{aligned}
\varepsilon h_{-1} & =\id_{\Zset}\\
 h_{-1} \varepsilon + d_1 h_0 & =\id_{F_0}
\\
h_{i-1} d_i + d_{i+1} h_i & =\id_{F_i}, &i >0 .
\end{aligned}
\end{equation*}
The existence of the contracting homotopy for $F_\ast$ assures in particular the exactness of the resolution.
\end{defn}

Given a free resolution $F_\ast$, one can consider the chain complex of $\Zset$-modules (that is to say, Abelian groups) $C_\ast=(C_n,d_{C_n})_{n\in\Nset}$ defined by
$$
C_n=(\Zset \otimes_{\Zset G} F_\ast)_n, \quad n \geq 0
$$
(where $\Zset\equiv C_\ast(\Zset,0)$ is the chain complex with only one non-null $\Zset G$-module in dimension~$0$) with differential maps \mbox{$d_{C_n}: C_n \rightarrow C_{n-1}$} induced by \mbox{$d_n:F_n \rightarrow F_{n-1}$}.

Let us emphasize the difference between the chain complexes $F_\ast$ and $C_\ast=\Zset \otimes_{\Zset G} F_\ast$. The elements of $F_n$ ($n \geq 0$) can be seen as \emph{words} $\sum \lambda_i (g_i,z_i)$ where $\lambda_i \in \Zset$, $g_i \in G$ and~$z_i$ is a generator of $F_n$ (which is a \emph{free} $\Zset G$-module). On the other hand, the associated chain complex $C_\ast=\Zset \otimes_{\Zset G} F_\ast$ of Abelian groups has elements in degree $n$ of the form $\sum \lambda_i z_i$
with $\lambda_i \in \Zset$ and $z_i$  a generator of the free $\Zset$-module $C_n$.

Although the chain complex of $\Zset G$-modules $F_\ast$ is acyclic, $C_\ast=\Zset \otimes_{\Zset G} F_\ast$ is in general not exact and  its homology groups are thus not null. An important result in homology of groups claims that these homology groups are independent of the chosen resolution for~$G$.

\begin{thm}(\cite{Bro82})
Let $G$ be a group and $F_\ast$, $F'_\ast$ two free resolutions of $G$. Then
$$
H_n(\Zset \otimes_{\Zset G} F_\ast) \cong H_n(\Zset \otimes_{\Zset G} F'_\ast) \quad \mbox{for all } n \in \Nset.
$$
\end{thm}

The  hypothesis that $F_\ast$ and $F'_\ast$ are free can in fact be relaxed; it suffices for the modules~$F_\ast$ and $F'_\ast$ to be \emph{projective}.
This theorem leads to the following definition.

\begin{defn}
Given a group $G$, the \emph{homology groups} $H_n(G)$ are defined as
$$H_n(G)=H_n(\Zset \otimes_{\Zset G} F_\ast), \quad n \in \Nset$$
where $F_\ast$ is any free (or projective) resolution for $G$.
\end{defn}

Let $G$ be a group, how can we determine a free
resolution $F_\ast$? One approach is to consider the \emph{Bar resolution} $B_\ast=\barc_\ast(G)$ (explained, for instance, in \cite{Mac63}) whose associated chain complex $\Zset \otimes_{\Zset G} B_\ast$ can be viewed as the chain complex of the Eilenberg-MacLane space $K(G,1)$ (see~\cite{Bro82} for details).
The homology groups of $K(G,1)$ are those of the group $G$ and this space has a big structural richness. But it has a serious drawback: its size. If $n>1$, then $K(G,1)_n=G^n$. In particular, if $G=\Zset$, the space $K(G,1)$ is of infinite type in each dimension. This fact is an important obstacle to
 using $K(G,1)$ as a means for computing  the homology groups of $G$. It would be therefore convenient to  construct \emph{smaller} resolutions.

For some particular cases, small (or minimal) resolutions can be directly constructed. For instance, let $G$ be the cyclic group of order $m$ with generator~$t$, $G=C_m$. The resolution~$F_\ast$ for~$G$
$$
\cdots \stackrel{t-1}{\longrightarrow} \Zset G \stackrel{N}{\longrightarrow} \Zset G  \stackrel{t-1}{\longrightarrow} \Zset G   {\longrightarrow} \Zset \longrightarrow 0
$$
where $N$ denotes the \emph{norm element} $1+ t + \cdots + t^{m-1}$ of $\Zset G$,
produces the chain complex of Abelian groups
$$
\cdots \stackrel{0}{\longrightarrow} \Zset  \stackrel{m}{\longrightarrow} \Zset   \stackrel{0}{\longrightarrow} \Zset     \longrightarrow 0
$$
and therefore
$$
H_i(G)=\left\{
\begin{array}{ll}
	\Zset & \mbox{if } \ i=0 \\
	\Zset / m\Zset & \mbox{if } \ i \mbox{ is odd} \\
	0 & \mbox{if } \ i \mbox{ is even and }i>0.
\end{array}
\right.
$$

But in general it is not so easy to obtain a resolution for a group~$G$, and in fact this problem provides an interesting research field where many papers and works have appeared trying to determine resolutions for different kinds of groups. As we will see later, the GAP package HAP has been designed
as a tool for constructing resolutions for a wide variety of groups. On the other hand, our work shows that the \emph{effective homology} method, introduced in the following section, could also be helpful in order to compute the homology of some groups.

\subsection{Effective homology}
\label{sec:efhm}

We now present the general ideas of the
effective homology method, devoted to the computation of homology groups of \emph{spaces}. See~\cite{RS02} and
\cite{RS06} for more details.

\begin{defn}
\label{def:red} A \emph{reduction} $\rho$ between two chain complexes
\mbox{$C_\ast=(C_n,d_{C_n})_{n \in \Nset}$} and
\mbox{$D_\ast=(D_n,d_{D_n})_{n \in \Nset}$} (which is denoted \mbox{$\rho: C_\ast \rrdc
D_\ast$}) is a triple $(f,g,h)$
where: (a)~the components $f$ and $g$ are chain complex morphisms
$f: C_\ast \rightarrow D_\ast$ and $g: D_\ast \rightarrow C_\ast$;
(b)~the component $h$ is a homotopy operator $h:C_\ast\rightarrow
C_{\ast+1}$ (a graded group morphism of degree~{+1}); (c)~the
following relations are satisfied:
$f  g = \id_D$;  $ d_C h + h  d_C = \id_C - g f $;
{$f  h = 0$;} $h   g = 0$; $h   h = 0$.
\end{defn}

These properties express that $C_\ast$ is the direct sum of
$D_\ast$ and an acyclic complex. This decomposition
is simply $C_\ast=\kr  f \oplus \im  g$, with $\im g\cong D_\ast$
and $H_\ast(\kr  f )=0$. In particular, this implies that the
graded homology groups \(H_\ast(C_\ast)\) and \(H_\ast(D_\ast)\)
are canonically isomorphic.

\begin{rmk}
A reduction is in fact a particular case of chain equivalence in the classical sense (see~\cite{Mac63}, page 40), where the homotopy operator on the small chain complex $D_\ast$ is the null map.
\end{rmk}

\begin{defn}
A \emph{(strong chain) equivalence} $\varepsilon$ between two
chain complexes $C_\ast$ and $D_\ast$, denoted by $\varepsilon:
C_\ast \lrrdc D_\ast$, is a triple $(B_\ast,\rho_1,\rho_2)$ where
$B_\ast$ is a chain complex, and $\rho_1$ and $\rho_2$ are
reductions $\rho_1:B_\ast \rrdc C_\ast$ and $\rho_2:B_\ast\rrdc D_\ast$.
\end{defn}

\begin{rmk}
We need the notion of \emph{effective} chain complex: it is essentially a free chain complex
$C_\ast$ where each group $C_n$ is finitely generated, and a
provided algorithm returns a (distinguished) $\Zset$-basis in each
degree $n$; in particular, its homology groups are elementarily computable (for details, see~\cite{RS02}).
\end{rmk}

\begin{defn}
An \emph{object with effective homology} $X$ is a quadruple $(X,C_\ast(X), HC_\ast, \varepsilon)$ where
$C_\ast(X)$ is a chain complex canonically associated with $X$ (which allows us to study the homological nature of $X$), $HC_\ast$ is an effective chain complex, and $\varepsilon$ is an equivalence $\varepsilon: C_\ast (X) \lrrdc HC_\ast$.
\end{defn}

It is important to understand that in general the \(HC_\ast\)
component of an object with effective homology is \emph{not} made
of the homology groups of \(X\); this component \(HC_\ast\) is a
free \(\Zset\)-chain complex of finite type, in general with a
non-null differential, whose homology groups $H_\ast(HC_\ast)$ can be determined by
means of an elementary algorithm. From the equivalence $\varepsilon$
one can deduce the isomorphism $H_\ast(X):=H_\ast(C_\ast(X))\cong H_\ast(HC_\ast)$, which allows one to compute the homology groups of the initial space \(X\).
In this way, the notion of object with effective homology provides a method to compute homology groups of complicated spaces by means
of homology groups of effective complexes.

The effective homology technique is based
on the following idea: given some topological spaces $X_1, \ldots,
X_n$, a topological constructor $\Phi$ produces a new topological
space~$X$. If effective homology versions of the spaces $X_1,
\ldots, X_n$ are known, then one should be able to build an effective homology
version of the space $X$, and this version would allow us to
compute the homology groups of $X$. A typical example of this kind of situation is the loop space
constructor. Given a {$1$-reduced} simplicial set $X$ with effective
homology, it is possible to determine the effective homology of
the loop space $\Omega(X)$, which in particular allows one to
compute the homology groups $H_\ast(\Omega(X))$. Moreover, if $X$
is \mbox{$m$-reduced}, this process may be iterated $m$ times, producing
an effective homology version of $\Omega^k(X)$, for $k \leq m$.
Effective homology versions are also
known for classifying spaces or total spaces of fibrations, see
\cite{RS06} for more information.

All these constructions have been implemented in the Kenzo system (\cite{Kenzo}), a Common Lisp program which makes use of the effective homology method to determine homology groups of complicated spaces; it has obtained some results (for example homology groups of iterated loop spaces of a loop space modified by a cell attachment, components of complex Postnikov towers, etc.) which had never been determined before. Furthermore, Kenzo implements Eilenberg-MacLane spaces $K(G,n)$ for every $n$ but only for $G=\Zset$ and $G=\Zset/ 2 \Zset$ (these spaces appear  in different constructions of Algebraic Topology), although in principle it is not designed to determine the homology of groups and it does not know how to work with resolutions.

These ideas suggest that the effective homology technique and the Kenzo program
should have a role in the computation
of the homology of a group $G$. To this end, we can consider the Eilenberg-MacLane space $K(G,1)$, whose homology groups coincide with those of $G$. The size of this space makes it difficult to calculate the groups in a direct way, but it is possible to operate with this simplicial set making use of the \emph{effective homology} technique: if we construct the effective homology of $K(G,1)$ then we would be able to \emph{efficiently} compute the homology groups of $K(G,1)$, which are those of $G$.
Furthermore, it should be possible to extend many group theoretic constructions to  effective homology  constructions of Eilenberg-MacLane spaces. We thus introduce the following definition.

\begin{defn}
A group $G$ is a \emph{group with effective homology} if $K(G, 1)$ is a simplicial set with
effective homology.
\end{defn}

The problem is, given a group $G$, how can we determine the effective homology of $K(G,1)$? If the group $G$ is finite, the simplicial set $K(G,1)$ is effective too, so that it can be considered with effective homology in a trivial way. However, the enormous size of this space makes it difficult to obtain real calculations, and therefore even in this case we will try to obtain an equivalence $C_\ast(K(G,1))\lrrdc E_\ast$ where $E_\ast$ is effective and (much) smaller than the initial complex. Section~\ref{sec:algorithm} of this paper presents an algorithm that computes this desired equivalence provided that the group $G$ is endowed with a finite resolution.

\section{Effective homology of a group from a resolution}
\label{sec:algorithm}

This section is devoted to an algorithm computing the effective homology of a group~$G$ given a (small) free $\Zset G$-resolution. This algorithm was the main theoretical result included in the work \cite{RER09}.

\begin{alg} %
\label{alg}
\emph{Input:} a group $G$ and a free (augmented) finite type resolution $F_\ast$ for $G$ with a contracting homotopy $h$. \\
 \emph{Output:} the effective homology of the space $K(G,1)$, that is to say, a (strong chain) equivalence $C_\ast(K(G,1)) \lrrdc E_\ast$ where $E_\ast$ is an effective chain complex.
\end{alg}

Algorithm~\ref{alg} has been implemented in Common Lisp enhancing the Kenzo system. We will see some examples of use of these new programs in Section \ref{sec:Kenzo}. A brief description of the construction of the algorithm is included in the following paragraphs. For more details, see \cite{RER09}.

\begin{proof}
We begin by considering the Bar resolution  $B_\ast=\barc_\ast(G)$ for $G$,  with augmentation~$\varepsilon'$ and contracting homotopy $h'$ (the definition of these maps can be found in~\cite{Bro82}). As $B_\ast$ and the given resolution $F_\ast$ are free resolutions for $G$, it is well known (see ~\cite{Bro82}) that one can explicitly construct morphisms of chain complexes of \mbox{$\Zset G$-modules} \mbox{$f: B_\ast \rightarrow F_\ast$} and $g: F_\ast \rightarrow B_\ast$ which are homotopy equivalences. Moreover, one can construct graded morphisms of $\Zset G$-modules
$$k: F_\ast \rightarrow F_{\ast+1}, \quad \quad k': B_\ast \rightarrow B_{\ast+1}$$
such that $d_Fk + kd_F=\id_F - fg$ and $d_Bk'+k'd_B=\id_{B} - gf$.

We have therefore a homotopy equivalence (in the classical sense):
$$
\xymatrix {
B_\ast \ar @(l,u) []^{k'} \ar @/^/[rr]^f & & F_\ast \ar @/^/ [ll]^g \ar @(r,u) []_{k}
}
$$
in which the four components $f$, $g$, $k$ and $k'$ are compatible with the action of the group~$G$.

If we now apply the functor $\Zset \otimes_{\Zset G}-$, which is additive, we obtain an equivalence of chain complexes (of $\Zset$-modules):
$$
\xymatrix {
\Zset \otimes_{\Zset G}B_\ast \ar @(l,u) []^{k'} \ar @/^/[rr]^f & & \Zset \otimes_{\Zset G}F_\ast \ar @/^/ [ll]^g \ar @(r,u) []_{k}
}
$$
where both chain complexes provide us the homology of the initial group $G$, that is,
$$
H_\ast(\Zset \otimes_{\Zset G}B_\ast)\cong H_\ast(\Zset \otimes_{\Zset G}F_\ast)\equiv H_\ast(G).
$$

In order to obtain a strong chain equivalence (in other words, a pair of reductions, following the framework of effective homology), we make use of the mapping cylinder construction (see \cite{Wei94}). This allows one to produce a (strong chain) equivalence
$$
\Zset \otimes_{\Zset G}B_\ast \stackrel{\rho'}{\lrdc} \cil(f)_\ast \stackrel{\rho}{\rrdc} \Zset \otimes_{\Zset G}F_\ast.
$$
The definitions of the different components of both reductions are included in \cite{RER09}.

Now we recall that the left chain complex \mbox{$\Zset \otimes_{\Zset G}B_\ast$} is equal to $C_\ast(K(G,1))$. On the other hand, if we suppose that the initial resolution $F_\ast$ is of finite type (and small), then the right chain complex $\Zset \otimes_{\Zset G}F_\ast\equiv E_\ast$ is effective (and small too), so that we have obtained the desired
effective homology of $K(G,1)$,
$$C_\ast(K(G,1)) \lrrdc E_\ast.$$
\end{proof}

This strong chain equivalence makes it possible to determine the homology groups of~$G$, and, what is more useful, once we have $K(G,1)$ with its effective homology we could apply different constructors and obtain the effective homology of the results. This could allow one, for instance, to determine the homology of some groups (obtained from other inicial groups with effective homology) without constructing a resolution for them. Some fields of application of our algorithm are introduced in the following section.

\section{Applications}
\label{sec:applications}

\subsection{$2$-types}

Let $A$ be an Abelian group and $G$ a group acting on $A$; a \emph{$2$-type} for $G$ and $A$ is a (topological) space with $\pi_1(X)=G$, $\pi_2(X)=A$, and $\pi_n(X)=0$ for all $n \geq 3$; the computation of the homology groups of these spaces is a difficult problem in the field of group homology (\cite{Ell92}). It is well known that a $2$-type $X$ for $G$ and $A$ corresponds to a cohomology class $[f]$ in $H^3(G,A)$, and there exists a fibration
$$
K(A,2) \hookrightarrow X \rightarrow K(G,1).
$$

The theoretical existence of this fibration can be made \emph{constructive} as follows, when the action of $G$ on $A$ is trivial. A cohomology class $[f]$ is given by a $3$-cocycle $f$, which is a map $f: K(G,1)_3 \rightarrow A$ (satisfying some properties). This map induces a simplicial morphism $f:K(G,1)\rightarrow K(A,3)$, which can be composed with the universal fibration (see \cite{May67}) $K(A,2) \hookrightarrow E \rightarrow K(A,3)$ in order to construct the desired fibration. In this way, we obtain a twisting operator $\tau_f: K(G,1)_\ast \rightarrow K(A,2)_{\ast-1}$ which allows one to express the total space $X$ as a twisted Cartesian product
$$ X= K(A,2) \times_f K(G,1).$$


Supposing now that the group $G$ is given with a finite resolution, our Algorithm \ref{alg} can be applied in order to produce the effective homology of $K(G,1)$. Analogously, provided a finite resolution for $A$, we can determine the effective homology of $K(A,1)$. Since $K(A,1)$ is a simplicial Abelian group one can apply the classifying space constructor $\W$ that gives us $\W(K(A,1))=K(A,2)$, which is also a simplicial Abelian group. Furthermore, the  effective homology version of the classifying space constructor $\W$ (see \cite{RS06} for details) provides us the effective homology of the space $K(A,2)$ from the effective homology of $K(A,1)$ (iterating the process, $K(A,n)=\W(K(A,n-1))$ has effective homology for every $n \in \Nset$). In this way, both spaces $K(A,2)$ and $K(G,1)$ are \emph{objects with effective homology}. Finally, the effective homology version for a fibration (described also in \cite{RS06}), makes use of the effective homologies of $K(A,2)$ and $K(G,1)$ and of the twisting operator $\tau_f: K(G,1)_\ast \rightarrow K(A,2)_{\ast-1}$ and gives us the effective homology of the total space $X= K(A,2) \times_f K(G,1)$. In particular, this leads to the desired homology groups of the $2$-type $X$.

\begin{alg} %
\label{alg-2types}
 \emph{Input:}
\begin{itemize}
\item    an Abelian group $A$ with a free resolution $F_{A_\ast}$ of finite type (with a contracting homotopy);
\item    a group $G$ (acting trivially on $A$) with a free resolution $F_{G_\ast}$ of finite type (with a contracting homotopy);
\item    a cohomology class $[f]\in H^3(G,A)$ given by a $3$-cocycle $f: K(G,1)_3\rightarrow A$.
\end{itemize}
\emph{Output:} the effective homology of the $2$-type associated with $f$, $X= K(A,2) \times_f K(G,1)$.
\end{alg}

This algorithm have been implemented in Common Lisp as part of our new module for the Kenzo system dealing with homology of groups.
See Section \ref{kenzo:2types} for some examples of calculations.

If the group $G$ acts non-trivially on $A$, an action $K(G,0) \times K(A,2) \rightarrow K(A,2)$ must also be considered in the fibration $K(A,2) \hookrightarrow X \rightarrow K(G,1)$. The explicit construction of the twisting operator which describes the fibration cannot be obtained as easily as in the previous case, and a more deep study of the fibration is necessary. It should be done as a further work.

\subsection{Central extensions}

Let $0 \rightarrow A \rightarrow E \rightarrow G \rightarrow 1$ be a central extension of groups (that is, $A$ is an Abelian group and $G$ acts on
$A$ in a trivial way). Then, it is well-known (see \cite{Bro82}) there exists a set-theoretic map $f: G \times G \rightarrow A$ which satisfies:
\begin{enumerate}
\item $f(g, 1) = 0 =f(1, g)$
\item $f(gh, k) = f(h, k) -f(g, h) + f(g, hk)$
\end{enumerate}

In addition, the initial extension is equivalent to another extension
$$0 \rightarrow A \rightarrow A \times_fG \rightarrow G \rightarrow 1$$
where the elements of $A \times_f G$ are pairs $(a,g)$ with $a \in A$ and $g \in G$, and the group law is defined by
$$(a_1,g_1)(a_2,g_2) \equiv (a_1 + a_2 +f(g_1,g_2),g_1g_2).$$
The set-theoretic map $f$ is called the \emph{$2$-cocycle} of the extension, since it corresponds to a map $f:K(G,1)_2 \rightarrow A$ in $H^2(G,A)$.

Very frequently, the groups $G$ and $A$ are not complicated and their homology groups are known. On the contrary, the homology groups of $E \cong A \times_fG$ are not always easy to obtain. The effective homology technique and our Algorithm \ref{alg-centralext} will provide a method computing the desired homology groups of $E$ from finite resolutions of $G$ and $A$. In this way, it will not be necessary to determine a finite resolution for $E$.

As explained in a previous work of the second author of this paper (see \cite{Rub97}), given a $2$-cocycle $f$ defining a central extension of a group $G$ by an Abelian group $A$, one can (explicitly) construct a fibration
 $$K(A,1) \hookrightarrow X \rightarrow K(G,1)$$
 where the total space $X$ can be seen as a twisted Cartesian product $K(A, 1) \times_\tau K(G, 1)$. Furthermore, it can be proved that this space is in fact isomorphic to the Eilenberg-MacLane space $K(A \times_fG, 1)$, whose homology groups are those of the group $A \times_fG \cong E$. The simplicial morphisms $\Phi: K(A\times_f G,1) \rightarrow K(A,1) \times_\tau K(G,1)$ and $\Phi^{-1}: K(A,1) \times_\tau K(G,1) \rightarrow K(A\times_f G,1)$ can be found in \cite{Rub97}.

On the other hand, in the case where both the fiber and base spaces of the fibration, $K(A,1)$ and $K(G,1)$, are objects with effective homology, the effective homology version of a fibration (see \cite{RS06}) provides the effective homology of the total space $K(A,1) \times_\tau K(G, 1)$, which in particular will make it possible to obtain the homology groups of $E$. Finally, if the groups $G$ and $A$ are given with finite (small) resolutions, our Algorithm \ref{alg} provides the necessary effective homologies of $K(G,1)$ and $K(A,1)$. We obtain therefore the following algorithm.

\begin{alg} 
\label{alg-centralext}
\emph{Input:} \begin{itemize}
\item groups $G$ and $A$ ($A$ is Abelian) and for both of them free resolutions of finite type with the corresponding contracting homotopies (or more generally, $G$ and $A$ with effective homology);
\item a $2$-cocycle $f$ defining an extension of $G$ by $A$.
\end{itemize}
\emph{Output:} the effective homology of the central extension $A \times_f G$.
\end{alg}

This algorithm has also been implemented in Common Lisp and in particular it allows us to determine the homology groups of central extensions of finitely generated Abelian groups. In Section \ref{kenzo:centralext} we include some examples of calculations.

\section{New modules for Kenzo and experimental results}
\label{sec:Kenzo}

As already mentioned in Section \ref{sec:efhm}, Kenzo~(\cite{Kenzo}) is a Common Lisp program devoted to Symbolic Computation in Algebraic Topology, developed by Francis Sergeraert and some co-workers. This system makes use of the effective homology method to determine homology groups of complicated spaces and has obtained some results which had never been determined before. In principle Kenzo was not intended to compute homology of groups but we have enhanced this system with a new module dealing with groups, resolutions, and Eilenberg-MacLane spaces (which where already implemented in Kenzo for the particular cases $\Zset$ and $\Zset / 2 \Zset$) and we have written the corresponding programs implementing Algorithm \ref{alg}, which produces the effective homology of the space~$K(G,1)$ given a finite resolution for the group $G$. Since the construction of a finite resolution for a group is not always an easy task, we have allowed Kenzo to connect with the GAP package HAP and obtain a resolution from it. Furthermore, as already announced, we provide programs which determine the homology groups of some $2$-types and central extensions.

\subsection{Interoperating with GAP}

GAP~(\cite{GAP}) is a system for computational discrete algebra with particular emphasis on Computational Group Theory. In our work we consider the HAP homological algebra library (\cite{HAP}) for use with GAP; it was written by Graham Ellis and is still under development.
The initial focus of HAP is on computations related to the cohomology of groups. A range of finite and infinite groups are handled, with  particular emphasis on integral coefficients. It also contains some functions for the integral (co)homology of: Lie rings, Leibniz rings, cat-1-groups and digital topological spaces. And in particular, HAP allows one to obtain (small) resolutions of many different groups, although it does not implement the Bar resolution nor Eilenberg-MacLane spaces~$K(G,1)$.

As presented in \cite{RER09}, in a joint work with Graham Ellis, we have developed a new module making it possible to export resolutions from HAP and import them into Kenzo. As interchange language we have used OpenMath~(\cite{OM}), an XML standard for representing mathematical objects. There exist OpenMath translators from several Computer Algebra systems, and in particular GAP includes a package~(\cite{GAPOM}) which produces OpenMath code from some GAP elements (lists, groups...). We have extended this package in order to represent resolutions, including a new GAP function which provides the OpenMath code of these elements. In \cite{RER09} a detailed description of our OpenMath representation for resolutions can be found.

The communication between HAP and Kenzo is done as follows: given a group $G$, the system HAP produces a $\Zset G$-resolution (including the homotopy operator). This resolution can be automatically translated to OpenMath code thanks to our new function added to the OpenMath package for GAP, and this code is written in a text file. Then Kenzo imports the file (and translates the OpenMath code into Kenzo elements thanks to the corresponding parser) so that one can use the resolution directly in Kenzo without the need of programming it in Common Lisp. Once the resolution is defined in Kenzo, we can use it to determine the effective homology of $K(G,1)$ as explained in Section~\ref{sec:algorithm}. In this way, if the construction of a resolution for a group $G$ is complicated, we can avoid  programming it by hand; it will be automatically implemented in Kenzo by obtaining it from HAP.

\subsection{Computations with K(G,n)'s}
\label{sec:kgns}

Let $G$ be a group and let us suppose that Kenzo knows a resolution for it (for some particular groups the system can construct it directly; for others, it could obtain it from HAP). Making use of our main Algorithm \ref{alg} one can determine the effective homology of the simplicial Abelian group $K(G,1)$, and in particular, compute its homology groups.

Let us consider, for instance, $G=C_5$, the cyclic group of order $5$. As already seen in Section \ref{sec:group-hmlg}, in this case it
is not difficult to construct a small resolution $F_\ast$ of $G$. The group can be built in Kenzo with the function \texttt{cyclicGroup}; the program computes automatically the well-known resolution for $G$ (coded as a reduction $F_\ast \rrdc \Zset$) and stores it in the slot \texttt{resolution} of the group.

\vskip0.2cm
\scriptsize
\begin{verbatim}
> (setf C5 (cyclicGroup 5))
[K1 Abelian-Group]
> (resolution C5)
[K10 Reduction K2 => K5]
\end{verbatim}

\vskip0.2cm
\normalsize
This resolution is then used by our programs, following Algorithm \ref{alg}, to determine the effective homology of $K(G,1)$. The corresponding homotopy equivalence is spontaneously computed and stored in the slot \texttt{efhm}. In this way one can obtain the homology groups of this Eilenberg-MacLane space.

\vskip0.2cm
\scriptsize
\begin{verbatim}
> (setf KC51 (K-G-1 C5))
[K11 Abelian-Simplicial-Group]
> (efhm KC51)
[K50 Homotopy-Equivalence K11 <= K40 => K31]
> (homology KC51 0 5)
Homology in dimension 0 :
Component Z
---done---
Homology in dimension 1 :
Component Z/5Z
---done---
Homology in dimension 2 :
---done---
Homology in dimension 3 :
Component Z/5Z
---done---
Homology in dimension 4 :
---done---
\end{verbatim}

\vskip0.2cm
\normalsize Moreover, since $G=C_5$ is Abelian, $K(G,1)$ is a simplicial Abelian group, and we can apply the classifying space constructor $\W$ (already implemented in Kenzo) which gives us $\W(K(G,1))=K(G,2)$, a new simplicial Abelian group \emph{with effective homology}.

\vskip0.2cm
\scriptsize
\begin{verbatim}
> (setf KC52 (classifying-space KC51))
[K51 Abelian-Simplicial-Group]
> (efhm KC52)
[K190 Homotopy-Equivalence K51 <= K180 => K176]
> (homology KC52 3 6)
Homology in dimension 3 :
---done---
Homology in dimension 4 :
Component Z/5Z
---done---
Homology in dimension 5 :
---done---
\end{verbatim}

\vskip0.2cm
\normalsize
Iterating the process, $K(G,n)=\W(K(G,n-1))$ has effective homology for every $n \in \Nset$. Our new Kenzo function \texttt{K-Cm-n} allows us to directly construct $K(C_m,n)$; we observe that the slot \texttt{efhm} is automatically created.

\vskip0.2cm
\scriptsize
\begin{verbatim}
> (setf KC42 (K-Cm-n 4 2))
[K204 Abelian-Simplicial-Group]
> (efhm KC42)
[K378 Homotopy-Equivalence K204 <= K368 => K364]
> (homology KC42 4)
Homology in dimension 4 :
Component Z/8Z
---done---
\end{verbatim}

\vskip0.2cm
\normalsize The construction of Eilenberg-MacLane spaces $K(G,n)$ for every cyclic group $G=C_m$ (with the corresponding effective homology) is an important enhancing of the Kenzo system, which previously was only able to deal with cases $G=\Zset$ and $G=\Zset/2 \Zset=C_2$. The homology groups obtained for some of these new spaces have been tested comparing them with the results shown in Alain Cl\'ement's thesis (\cite{Cle02}). It is important to stress that Cl\'ement's tables, computed by using a direct algorithm created by Henry Cartan (see \cite{Car54}), contain much more groups than those that can be computed with our programs in its current state. Nevertheless, Cl\'ement's tables give only the homology groups of the spaces, while our approach provides the \emph{effective homology}. Our information is much more complete, giving access to geometrical generators of the homology and, in fact, fully solving the \emph{homological problem} of these groups (see \cite{RS06}). And, perhaps more important, our programs allow us to continue working with the corresponding $K(G,n)$, to produce new interesting topological spaces and to determine their homology groups. The information computed by Cl\'ement is not enough to carry out this further work.

The same technique explained for cyclic groups can  be used to compute the effective homology of
spaces $K(G,n)$, where $G$ is a finitely generated Abelian group.
In this case, the homology of $K(G,n)$ is one of the main ingredients
to compute homotopy groups of spaces (see~\cite{RS02} and~\cite{RS06} for details).

\subsection{An example of homology of a $2$-type}
\label{kenzo:2types}
Let us consider now $G=C_3$ the cyclic group of order~$3$. Let $A=\Zset/3\Zset$ be the Abelian group of three elements with trivial $G$-action (the groups $G$ and $A$ are in fact isomorphic; different notations are used to distinguish multiplicative and additive operations). Then the third cohomology group of $G$ with
coefficients in $A$ is
$$
H^3(G,A) = \Zset/3\Zset.
$$

The classes $[f]$ of this cohomology group correspond to $2$-types with $\pi_1=G$,
$\pi_2=A$, and one such $2$-type $X$ can be seen as a twisted Cartesian product  $ X= K(A,2) \times_f K(G,1)$. It can be constructed by Kenzo in the following way:

\vskip0.2cm
\scriptsize
\begin{verbatim}
> (setf KC31 (K-Cm-n 3 1))
[K380 Abelian-Simplicial-Group]
> (setf chml-clss (chml-clss KC31 3))
[K427 Cohomology-Class on K407 of degree 3]
> (setf tau (zp-whitehead 3 KC31 chml-clss))
[K442 Fibration K380 -> K428]
> (setf X (fibration-total tau))
[K448 Kan-Simplicial-Set]
\end{verbatim}

\vskip0.2cm
\normalsize
As seen in the previous section, $K(A,2)$ and $K(G,1)$ are objects with effective homology. From the two equivalences $C_\ast(K(A,2)) \lrrdc E_\ast$ and $C_\ast(K(G,1)) \lrrdc E'_\ast$, Kenzo knows how to construct the effective homology of the twisted Cartesian product $X = K(A,2) \times_f K(G,1)$, which allows one to determine its homology groups.

\vskip0.2cm
\scriptsize
\begin{verbatim}
> (efhm X)
[K660 Homotopy-Equivalence K448 <= K650 => K646]
> (homology X 5)
Homology in dimension 5 :
component Z/3Z
---done---
\end{verbatim}

\vskip0.2cm
\normalsize
In the same way, the homology groups of $X = K(A,2) \times_f K(G,1)$ can be determined for all groups $A$ and $G$ with given (small) resolutions and cohomology classes~$[f]$ in~$H^3(G,A)$. 
Up to now, only the homology of finitely presented groups has been considered, restricting the kind of $2$-types that can be studied with our methods, since only spaces with Abelian fundamental group would be in its scope. The range of groups which can be considered is considerably enlarged with the central extension constructions, as explained in the following subsection.

\subsection{Central extensions}
\label{kenzo:centralext}

Let us introduce an interesting example of central extension extracted from \cite{Lea91}. Let $E$ be the group defined by the following presentation:
$$E= <x,y, z | x^p = y^p = z^{p^{n-2}} = [x, z] = [y,z] = 1; [x, y] = z^{p^{n-3}}>.$$

This group can be seen as a central extension of the groups
$$A = <z | z^{p^{n-2}} = 1),$$
isomorphic to the cyclic group with $p^{n-2}$ elements, and
$$G = <x, y| x^p=y^p=[x, y]= 1>,$$
which is the direct sum of two cyclic groups of cardinality $p$.
A $2$-cocycle of the extension is defined by
$$f(x^{p_1}y^{q_1} ,x^{p_2}y^{q_2} ) = z^{q_1p_2(p-1)p^{n-3}}.$$

As already explained, the group $A\cong C_{p^{n-2}}$ has effective homology. On the other hand, the effective homology of
$G\cong C_p \oplus C_p$ can be easily obtained from the effective homology of the cyclic group $C_p$ (a direct sum
of two groups can in fact be considered as a particular case of central extension, where the $2$-cocycle is trivial,
so that its effective homology can be computed given the effective homologies of the two factors).
In this way, Algorithm \ref{alg-centralext} can be applied
to obtain the effective homology of $E$ and then compute its homology groups.

Let us consider, for instance, $p=3$ and $n=4$. The following Kenzo instructions construct the group $E$.

\vskip0.2cm
\scriptsize
\begin{verbatim}
> (progn
     (setf p 3 n 4)
     (setf A (cyclicGroup (expt  p (- n 2))))
     (setf G (gr-crts-prdc (cyclicGroup p) (cyclicGroup p)))
     (setf cocycle #'(lambda (crpr1 crpr2)
                    (with-grcrpr (x1 y1) crpr1
                      (with-grcrpr (x2 y2) crpr2
                        (mod (* y1 x2 (1- p) (expt p (- n 3))) (expt p (- n 2)))))))
     (setf E (gr-cntr-extn A G cocycle)))
[K663 Group]
\end{verbatim}

\vskip0.2cm
\normalsize
The spaces $K(A,1)$ and $K(G,1)$ can be constructed with the function \texttt{K-G-1}; both of them are Abelian
simplicial groups with effective homology.

\vskip0.2cm
\scriptsize
\begin{verbatim}
> (setf KA1 (K-G-1 A))
[K664 Abelian-Simplicial-Group]
> (efhm KA1)
[K710 Homotopy-Equivalence K664 <= K700 => K691]
\end{verbatim}

\begin{verbatim}
> (setf KG1 (K-G-1 G))
[K711 Abelian-Simplicial-Group]
> (efhm KG1)
[K775 Homotopy-Equivalence K711 <= K765 => K745]
\end{verbatim}

\vskip0.2cm \normalsize
Given the effective homologies of $K(A,1)$ and $K(G,1)$, our Algorithm \ref{alg-centralext} returns the effective homology
of $K(E,1)$, which is then stored in the corresponding slot \texttt{efhm}.

\vskip0.2cm
\scriptsize
\begin{verbatim}
> (setf KE1 (K-G-1 E))
[K776 Simplicial-Group]
> (efhm KE1)
[K884 Homotopy-Equivalence K776 <= K870 => K866]
> (homology KE1 0 5)
Homology in dimension 0 :
Component Z
---done---
Homology in dimension 1 :
Component Z/3Z
Component Z/3Z
Component Z/3Z
---done---
Homology in dimension 2 :
Component Z/3Z
Component Z/3Z
---done---
Homology in dimension 3 :
Component Z/9Z
Component Z/3Z
Component Z/3Z
Component Z/3Z
---done---
Homology in dimension 4 :
Component Z/3Z
Component Z/3Z
Component Z/3Z
---done---
\end{verbatim}

\vskip0.2cm
\normalsize
In this way, one can determine the homology
groups of the central extension $E$. The computations obtained by our programs have been compared with Leary's theoretical
results for different values of $p$ and $n$; the same groups have been obtained by both methods. We can repeat here the discussion
made at the end of Subsection \ref{sec:kgns} with respect to Cl\'ement's computations for $H_\ast(K(G,n))$: Leary's methods give
more groups than our techniques, but with less information. In particular, our results allow us to compute the homology of $2$-types
whose fundamental groups are central extensions, while Leary's groups are not enough for this task.

\section{The inverse problem: recovering a resolution from the effective
homology of a group}

In Section \ref{sec:algorithm} we have presented an algorithm which, given a group $G$ with a free finite resolution $F_\ast$, constructs the effective homology of the simplicial Abelian group $K(G,1)$. This effective homology allows one to determine the homology groups of $G$ and, as seen in Sections \ref{sec:applications} and \ref{sec:Kenzo}, makes it possible to use the space $K(G,1)$ as initial data for some constructions in Algebraic Topology, computing in this way homology groups of other interesting objects.

We consider now the inverse problem: let $G$ be a group such that an equivalence
$$
C_\ast(K(G,1)) \lrrdc E_\ast
$$
is given, $E_\ast$ being a finite type chain complex of Abelian groups. Is it possible to obtain a finite free resolution for the group $G$? It seems, in principle, that the answer should be negative in the general case; since no condition is imposed on the arrows, they surely do not respect the $G$-action and, thus, it would not be possible to build a $\Zset G$-resolution. We have proved, however, that supposing that some additional conditions for the given chain equivalence are satisfied, one can construct the desired resolution with the corresponding contracting homotopy.

The algorithm we have developed makes use of the \emph{Basic Perturbation Lemma} (BPL), one of the fundamental results in Constructive Algebraic Topology. The general idea of this theorem is that given a reduction $\rho=(f,g,h):C_\ast \rrdc D_\ast$, if we modify the initial differential $d_C$ of the \emph{big} complex $C_\ast$ by adding some \emph{perturbation}, then it is possible to perturb the differential $d_D$ in the \emph{small} chain complex $D_\ast$ so that we obtain a new reduction between the perturbed complexes. But the result is not always true, a necessary condition must be satisfied: the composite function $h\circ d_C$ must be \emph{locally nilpotent}. An endomorphism $\alpha: C_\ast \rightarrow C_\ast$ is locally nilpotent if for every $x \in C_\ast$ there exists $m \in \Nset$ such that $\alpha^m(x)=0$. The condition of local nilpotency ensures the convergence of a formal series used in the BPL to build the perturbed differential on the small complex (and, in fact, to construct also all the arrows defining the new reduction between the perturbed complexes). The Basic Perturbation Lemma was discovered by Shih Weishu (\cite{Shi62}), and then generalized by Ronnie Brown in \cite{Bro67}. In its modern form it was formulated by Gughenheim (\cite{Gu72}) and its essential use in Kenzo has been documented in \cite{RS06}.

\begin{alg} 
\label{alg-inverse-problem}
\emph{Input:} \begin{itemize}
\item a group $G$;
\item a (strong) chain equivalence $C_\ast(K(G,1))\stackrel{\rho_1}{\lrdc} D_\ast \stackrel{\rho_2}{\rrdc} E_\ast$, where $\rho_1=(f_1,g_1,h_1)$, $\rho_2=(f_2,g_2,h_2)$, $E_\ast$ is an effective chain complex, and the composition $h_2g_1\partial_nf_1$ is locally nilpotent ($\partial_n$ is the face of index $n$ over the elements of $K(G,1)_n$, which can be extended to $C_n(K(G,1))$).
\end{itemize}
\emph{Output:} a free resolution $F_\ast$ for $G$ of finite type with a contracting homotopy $h$.
\end{alg}

\begin{proof}
Let us start by considering the universal fibration $K(G,0) \rightarrow K(G,0) \times_{\tau} K(G,1) \rightarrow K(G,1)$ (see \cite{May67}). The total space $K(G,0) \times_{\tau} K(G,1)$ is acyclic and one can construct a reduction
$$
C_\ast(K(G,0) \times_{\tau} K(G,1)) \rrdc \Zset
$$
where $\Zset$ represents the chain complex (of Abelian groups) $C_\ast(\Zset,0)$ with a unique non-null component $\Zset$ in dimension $0$.

On the other hand, one can consider the Eilenberg-Zilber theorem (\cite{EZ53}), which relates the chain complex of a Cartesian product with the tensor product of the chain complexes of the two components, and allows one to build a reduction
$$C_\ast(K(G,0) \times K(G,1)) \rrdc C_\ast(K(G,0)) \otimes C_\ast(K(G,1)).$$

Applying the BPL (it can be proved that the nilpotence condition is satisfied) we obtain a \emph{perturbed} reduction (this is in fact the \emph{twisted} Eilenberg-Zilber theorem, see \cite{May67})
$$C_\ast(K(G,0) \times_\tau K(G,1)) \rrdc C_\ast(K(G,0)) \otimes_t C_\ast(K(G,1))$$
where $C_\ast(K(G,0)) \otimes_t C_\ast(K(G,1))$ is a chain complex with the same underlying graded module as the tensor product $C_\ast(K(G,0)) \otimes C_\ast(K(G,1))$, but
its differential is modified to take account of the twisting operator $\tau$.

Now, from the given equivalence $C_\ast(K(G,1))\lrdc D_\ast \rrdc E_\ast$, it is not difficult to construct a new equivalence
$$
C_\ast(K(G,0)) \otimes C_\ast(K(G,1)) \lrdc C_\ast(K(G,0)) \otimes D_\ast \rrdc C_\ast(K(G,0)) \otimes E_\ast
$$
and applying again the BPL, provided that $h_2g_1{\partial}_nf_1$ is locally nilpotent, we obtain
 $$
C_\ast(K(G,0)) \otimes_t C_\ast(K(G,1)) \lrdc C_\ast(K(G,0)) \otimes_t D_\ast \rrdc C_\ast(K(G,0)) \otimes_t E_\ast.
$$

Finally, one can observe that $C_*(K(G,0)) \equiv \Zset G$ and composing the reductions \linebreak $
C_\ast(K(G,0) \times_{\tau} K(G,1)) \rrdc \Zset
$ and $C_\ast(K(G,0) \times_\tau K(G,1)) \rrdc C_\ast(K(G,0)) \otimes_t C_\ast(K(G,1))$ with the last equivalence,
we get a contracting
homotopy on $\Zset G \otimes_{t} E_*$ which is a resolution for $G$.
\end{proof}

As a first possible application of this theorem, one can consider the integer group \mbox{$G=\Zset$} and the well known effective homology of $K(\Zset,1)$, given by a reduction $C_\ast(K(\Zset,1))$ $\rrdc  C_\ast(S^1)$, where $S^1$ denotes a simplicial model for the sphere of dimension $1$. In this case it is not difficult to prove the desired condition, $h_2g_1\partial_nf_1$ is locally nilpotent, and therefore one can construct a finite resolution for $G=\Zset$, as a reduction $\Zset G \otimes_{t} C_*(S^1)\rrdc \Zset$.

A natural question which appears in this context is whether, given a group $G$ and an equivalence $C_\ast(K(G,1))\lrrdc E_\ast$ which has been obtained by means of our Algorithm \ref{alg} from a finite resolution $F_\ast$, the necessary condition of $h_2g_1\partial_nf_1$ being locally nilpotent is satisfied or not.
The answer is positive if the group $G$ and the resolution $F_\ast$ satisfy some particular properties. More concretely, we suppose that a \emph{norm} is defined on $G$ and it can be extended to $F_\ast$ in the following \emph{natural} way.

\begin{defn}
\label{defn:norm}
Let $G$ be a group. A \emph{norm} for $G$ is a map $||.||:G \rightarrow \Nset$ such that
\begin{itemize}
\item $||g||>0$ for each $g \in G$ and $||g||=0$ if and only if $g=1$;
\item $||g_1 g_2||\leq ||g_1||+||g_2||$ for all $g_1,g_2 \in G$.
\end{itemize}
We suppose that the resolution is reduced ($F_0=\Zset G$) and define $||.||: F_{0}=\Zset G\rightarrow \Nset $ as $||\sum \lambda_i g_i ||= \max \{||g_i||\}$. We say that the norm is \emph{compatible} with the resolution $F_\ast$ if for each $n \geq 1$ we can also define $||.||: F_n \rightarrow \Nset$ such that
\begin{itemize}
\item $||(g,z)||=||g||+ ||z||$ for all $g\in G$ and $z$ a generator of $F_n$;
\item there exists $i_n\in \Zset$ such that $||h_n(x)||\leq ||x||-i_n$ and $||d_{n+1}(x')||\leq ||x'||+i_n$ for all $x \in F_n, x'\in F_{n+1}$.
\end{itemize}
\end{defn}

The last condition introduces a control measure on the contracting homotopy $h$, with respect to the structure of the group, allowing us (as shown in the following result) to ensure in this case the convergence of the Basic Perturbation Lemma. Examples of resolutions with this kind of norm are the Bar resolution, the canonical small resolution for $G=\Zset$ and, for instance, the small resolutions for cyclic groups introduced in Section~\ref{sec:group-hmlg}.

\begin{thm}
\label{thm:norms}
Let $G$ be a group and $F_\ast$ a free resolution for $G$ with contracting homotopy~$h$. Let us suppose that $G$ is provided with a norm $||.||: G \rightarrow \Nset$ which is compatible with the resolution. Then the effective homology of $K(G,1)$ obtained from $F_\ast$ by our Algorithm~\ref{alg} satisfies the necessary condition of $h_2g_1\partial_nf_1$ being locally nilpotent, and therefore it is possible to construct a (new) free finite resolution for $G$.
\end{thm}

The proof of this theorem follows an inductive reasoning and is based on a deep study of the definition of the different components in the equivalence $C_\ast(K(G,1)) \lrrdc E_\ast$.

It is not difficult to observe that in this case the new resolution $F'_\ast$ given by Algorithm \ref{alg-inverse-problem} has the same structural components as the initial resolution $F_\ast$; in other words, $F_n=F'_n$ for all $n \in \Nset$. However, the differential and contracting homotopy maps could be different.

A final example of application of Algorithm \ref{alg-inverse-problem} and Theorem \ref{thm:norms} is the following one.

\begin{thm}
\label{thm:norms2}
Let $G,G'$ groups with free resolutions $F_\ast$ and $F'_\ast$ (with contracting homotopies $h$ and $h'$ respectively). Let us suppose that there exists norms on $G$ and $G'$ which are compatible with the corresponding resolutions. Then the effective homology of $G\oplus G'$ (obtained from those of $G$ and $G'$ as a particular case of central extension) satisfies that $h_2g_1\partial_nf_1$ is locally nilpotent, so that it is also possible to determine a resolution for the direct sum $G\oplus G'$.
\end{thm}

Again, we know the graded part of the output resolution, but it is still unknown if the differential and contracting homotopy constructed have some good geometrical behavior.

\section{Conclusions and further work}

In this paper we have defended this proposal: the \emph{geometric way}
for computing group homology can be sensible and fruitful. To this aim,
we have worked inside Sergeraert's \emph{effective homology}, and added
packages devoted to group homology in Sergeraert's Kenzo system.

In their current state our methods have a performance penalty
when compared with the more standard algebraic approach (based on
\emph{resolutions}).
Nevertheless, this claim is only true for computations reachable by previous means. Furthermore,
what is more important, to get available the homology of a group $G$ through
an Eilenberg-MacLane
space $K(G,1)$ with \emph{effective homology} allows us to use that
space for further
topological constructions. The poorer performance is therefore
balanced with the richer
information we get.

The paper illustrates our approach with concrete computer experiments
for general
Eilenberg-MacLane spaces $K(G,n)$, for central extensions of groups
and for $2$-types.
In the first two applications, the computer results have been compared
with previously
published works. In the case of the homology groups of $2$-types computed with
Kenzo, no comparison is possible, because no other source of results
is known by us.

Furthermore we have explored the problem of computing a resolution of $G$ from
the effective homology of $K(G,1)$, obtaining some partial algorithmic results
which have not yet been implemented.

Some of the lines opened in this paper have not been completely
closed, signaling
clear lines of further work. Starting from the end, the scope of the methods
to compute resolutions from effective homologies should be enlarged, and more
examples should be worked out. In particular, a comparison between the initial
resolutions and the ones constructed in the case of \emph{normed} groups should
be undertaken, trying to elucidate if our output resolution is better in some
geometrical sense.

In the area of $2$-types, the more important task would be to extend
our approach
to $2$-types with non-trivial action of the fundamental group. The main obstacle
here is to obtain a fibration expressed as a twisted Cartesian product, in order
to be able to apply the previous Kenzo infrastructure.

Another challenge consists in trying to get better algorithms from the efficiency point
of view, in such a way that our programs can compete with other
approaches. In particular,
we should improve the algorithm to construct the effective homology
from a resolution,
at least in certain cases, to obtain execution times closer to those of
the source system,
HAP. For finitely generated Abelian groups (which are the building
blocks to start
many of our constructions) a more direct approach, much more efficient,
could be extracted
from the original papers by Eilenberg and MacLane (\cite{EM53}, \cite{EM54}, and \cite{EM54b}).

Finally, the application of our methods for wider classes of groups
(for instance, extensions
beyond the central extensions dealt with in this paper) is likely possible and
surely an interesting research topic.

\section{Acknowledgements}

Thanks are due to Graham Ellis who collaborated with us in the connection of Kenzo and GAP
and in work dealing with $2$-types, as documented in \cite{RER09}.

\bibliographystyle{amsplain}
\bibliography{homology-of-groups-biblio}


\end{document}